\documentclass[11pt,oneside]{amsart}
\usepackage[T1]{fontenc}
\usepackage{lmodern}
\usepackage[hmargin=1.5in,vmargin=1.25in]{geometry}
\usepackage{etoolbox}
\newlength{\abstractmargin}
\makeatletter
\patchcmd{\abstract}{\leftmargin3pc}{\leftmargin\abstractmargin}{}{\PackageError{paper}{abstract patch failed}{}}
\makeatother
\makeatletter
\renewcommand{\@setaddresses}{%
  \par\addvspace{9pt}\begingroup\linespread{1}\footnotesize
  \def\author##1{}%
  \def\address##1##2{%
    \par\addvspace{3pt}\noindent
    \begin{minipage}{\linewidth}
      \raggedright\textbf{##1}\par\nobreak
      ##2
    \end{minipage}\par}%
  \addresses
  \endgroup}
\makeatother
\usepackage{amssymb,amsmath,mathtools}
\usepackage[expansion=false]{microtype}
\usepackage{needspace}
\usepackage{tikz}
\usepackage{xcolor}
\definecolor{linkblue}{RGB}{0,76,147}
\usepackage[colorlinks=true,allcolors=linkblue]{hyperref}
\hypersetup{
  pdftitle={Property F implies weak integrality for non-degenerate pivotal braided fusion categories},
  pdfauthor={},
  pdfsubject={Braided fusion categories and finite braid group images},
  pdfkeywords={property F, weak integrality, pivotal categories, pseudo-unitarity, Galois symmetry}
}

\newtheorem{theorem}{Theorem}[section]
\newtheorem{lemma}[theorem]{Lemma}
\newtheorem{proposition}[theorem]{Proposition}
\newtheorem{corollary}[theorem]{Corollary}
\newtheorem*{conjecture}{Conjecture}
\theoremstyle{remark}
\newtheorem{remark}[theorem]{Remark}

\newcommand{\C}{\mathcal C}
\newcommand{\M}{\mathcal M}
\newcommand{\B}{\mathcal B}
\newcommand{\E}{\mathcal E}
\newcommand{\D}{\mathcal D}
\newcommand{\Zc}{\mathcal Z}
\newcommand{\one}{\mathbf 1}
\newcommand{\FP}{\operatorname{FPdim}}
\newcommand{\Tr}{\operatorname{Tr}}
\newcommand{\tr}{\operatorname{tr}}
\newcommand{\ptr}{\operatorname{ptr}}
\newcommand{\Hom}{\operatorname{Hom}}
\newcommand{\End}{\operatorname{End}}
\newcommand{\Aut}{\operatorname{Aut}}
\newcommand{\sgn}{\operatorname{sgn}}
\newcommand{\id}{\mathrm{id}}
\newcommand{\Vect}{\mathrm{Vec}}

\newcommand{\Z}{\mathbb Z}
\newcommand{\Qb}{\overline{\mathbb Q}}
\newcommand{\ad}{\mathrm{ad}}
\newcommand{\rev}{\mathrm{rev}}
\newcommand{\cO}{\mathcal O}
\newcommand{\Gal}{\operatorname{Gal}}

\newcommand{\crossing}[2]{%
  \draw (#1+1,#2) .. controls (#1+1,#2+0.5) and (#1,#2+0.5) .. (#1,#2+1);
  \draw[white,line width=3.2pt] (#1,#2) .. controls (#1,#2+0.5) and (#1+1,#2+0.5) .. (#1+1,#2+1);
  \draw (#1,#2) .. controls (#1,#2+0.5) and (#1+1,#2+0.5) .. (#1+1,#2+1);}

\title[Property F implies weak integrality]{Property F implies weak integrality\\
for non-degenerate pivotal braided fusion categories}
\author{Justin Asher}
\address[Justin Asher]{Project Numina}
\author{Jonah Berean-Dutcher}
\address[Jonah Berean-Dutcher]{University of British Columbia, Vancouver, V6T 1Z1, Canada\newline
Institut des Hautes \'Etudes Scientifiques, 91440 Bures-sur-Yvette, France}
\author{Clay C\'ordova}
\address[Clay C\'ordova]{Leinweber Institute for Theoretical Physics \& Enrico Fermi Institute, University of Chicago}
\author{Julio Parra-Martinez}
\address[Julio Parra-Martinez]{Institut des Hautes \'Etudes Scientifiques, 91440 Bures-sur-Yvette, France}
\author{Richard Tuddenham}
\address[Richard Tuddenham]{Leinweber Institute for Theoretical Physics \& Enrico Fermi Institute, University of Chicago}
\hypersetup{pdfauthor={Justin Asher, Jonah Berean-Dutcher, Clay C\'ordova, Julio Parra-Martinez, Richard Tuddenham}}
\date{November 9, 2026}

\begin{document}
\pagestyle{plain}

\begin{abstract}
We prove that every non-degenerate pivotal braided fusion category with
property F is pseudo-unitary and weakly integral. In particular, every modular
category with property F is weakly integral. Neither sphericality nor
pseudo-unitarity is assumed. Passing to the Drinfeld centre extends the result
to degenerate categories whose dimensions are positive on the intersection of
the M\"uger centre with the adjoint subcategory; these include every braided
fusion category with a pseudo-unitary Galois conjugate.
Physically, the result connects the arithmetic of topological orders to
the computational power of anyon braiding: non-integer squared total
quantum dimension guarantees infinite braid group images, a necessary
condition for universal quantum computation.
\end{abstract}

\maketitle

\section{Introduction}

Let $\C$ be a braided fusion category over $\mathbb C$. For an object $X$ and
$n\ge2$, write $G_{X,n}$ for the image of
\[
\rho_{X,n}:B_n\longrightarrow\Aut(X^{\otimes n}),\qquad
\sigma_i\longmapsto\id_X^{\otimes(i-1)}\otimes c_{X,X}\otimes\id_X^{\otimes(n-i-1)}.
\]
The category has \emph{property F} if every $G_{X,n}$ is finite, and is
\emph{weakly integral} if $\FP(\C)\in\Z$.

\begin{conjecture}[Naidu--Rowell {\cite{NR}}]
A braided fusion category has property F if and only if it is weakly integral.
\end{conjecture}

Thus the conjecture relates finiteness of all braid group images to a single
arithmetic condition on the fusion category.
This is of particular physical interest for topological quantum
computation by anyon braiding \cite{K,FLWU}: the conjecture would provide
a simple arithmetic criterion identifying weakly integral
$(2+1)$-dimensional topological orders as incapable of universal
quantum computation by braiding anyons alone. Indeed, finite braid
group images cannot approximate arbitrary unitary gates
\cite{NSSFD,RW}.

The implication from weak integrality to property F is known for weakly
group-theoretical braided fusion categories \cite[Theorem~7.2]{GN}, extending the
group-theoretical case \cite{ERW}. Etingof, Nikshych and Ostrik asked whether there
exists a weakly integral fusion category that is not weakly group-theoretical
\cite[Question~2]{ENO2}; a negative answer would establish that implication in
general. In the other direction, infinite braid group images have been
established for several families of quantum-group categories \cite{FLW,LRW}.
We prove the implication from property F to weak integrality for non-degenerate
categories with a pivotal structure. Recall that $\C$ is \emph{pseudo-unitary} if $\dim(\C)=\FP(\C)$.

\begin{theorem}\label{thm:main}
Let $\C$ be a non-degenerate pivotal braided fusion category. If $\C$ has
property F, then $\C$ and all its Galois conjugates are pseudo-unitary, and
$\dim(\C)=\FP(\C)\in\Z$. In particular, every modular category with property F
is weakly integral.
\end{theorem}

Equivalently, a non-degenerate pivotal braided fusion category that is not
weakly integral admits a braid group representation with infinite image.

The Frobenius--Perron dimension is a character of the fusion ring. Since $\C$
is non-degenerate, it is the normalised $S$-matrix column of a simple object
$P$, and unitarity of the normalised $S$-matrix gives $\dim(\C)=|d_P|^2\FP(\C)$.
As $\dim(\C)\le\FP(\C)$, it suffices to show that $|d_P|\ge1$. The double
braiding of $P$ with itself has the largest trace that a character allows. We
show that finite braid images then force $|d_P|=\FP(P)$, which is at least one.

The main step is the trace argument in Section~\ref{sec:trace}. It first bounds
the categorical trace of every unitary in a braid algebra. The braid
$\sigma_1^2\sigma_2^2$ on three strands attains this bound; the equality case
makes $\sigma_1^2\sigma_2^{-2}$ act trivially wherever the trace sees it, and
computing its trace in two ways gives the required identity of dimensions.
Pivotal dimensions may be non-real, and the argument keeps track of their
phases.

Section~\ref{sec:prelim} records the transfer properties of property F and the
facts about pivotal dimensions used in the proof, and constructs the
Frobenius--Perron object. Section~\ref{sec:trace} establishes the trace bound
and the saturation lemma. The proof of Theorem~\ref{thm:main} is completed in
Section~\ref{sec:main}. Section~\ref{sec:degenerate} passes to the Drinfeld
centre. It extends the theorem to degenerate categories whose dimensions are positive on the
intersection of the M\"uger centre with the adjoint subcategory, which includes
every braided fusion category with a pseudo-unitary Galois conjugate. It then
derives a criterion for a single object to have infinite braid images.

\medskip
\noindent\textbf{Relation to recent work and transparency:}
We warmly acknowledge the recent preprint of Enyu Li and Yuze Ruan,
\emph{Property F and weak integrality of braided fusion categories} \cite{LR},
which establishes the principal result of this paper in comparable
generality. We congratulate the authors on their achievement and make
no claim of priority.

The argument presented here emerged from a human--AI collaboration
conducted independently of, and prior to our becoming aware of, \cite{LR}.
The investigation illustrates a broader approach to AI-assisted
mathematical research. Supported by a fellowship from \href{https://projectnumina.ai}{Project Numina},
it began with the Anyon Atlas \cite{AA}, a collection of modular data
generated by the authors, which provided a testing ground for autonomous
mathematical exploration. An initial round of agentic research produced
partial results that were subsequently understood, refined, and
strengthened by humans. Building on these improvements, a second round
of agentic research developed more general techniques and led to the
proof presented here. The resulting argument has also been formalized
in Lean.

Following the appearance of \cite{LR}, we have chosen to make our argument
publicly available in the interest of transparency, rather than to
claim credit for resolving the conjecture. The two proofs share the
same central idea: finite braid group images impose strong bounds on
categorical traces, and saturation of these bounds leads to the
equality of categorical and Frobenius--Perron dimensions. The arguments
differ in their implementation; the proof presented here uses an
explicit three-strand calculation to establish the required rigidity.

An accompanying \href{https://projectnumina.ai/article?slug=academic-loops}{blog post} documents the collaborative process,
including how human and AI contributions interacted in developing the
argument. The proof has also been \href{https://github.com/project-numina/mtc}{formalized in Lean}. We hope that sharing the proof,
its formal verification, and an account of their development will
offer a positive example of how human--AI collaboration can advance
scientific understanding.

\section{Preliminaries}\label{sec:prelim}

In the sequel, all fusion categories are over $\mathbb C$, and we suppress associativity and
unit constraints. We write $\one$ for the tensor unit, $\cO(\C)$ for the set of
simple objects up to isomorphism, and $N^W_Z=\dim\Hom(W,Z)$ for the
multiplicity of a simple object $W$ in an object $Z$. For a fusion subcategory
$\mathcal A$ of a braided category, $\mathcal A'$ denotes its centraliser; thus
$\C'$ is the M\"uger centre, and $\C$ is non-degenerate if $\C'=\Vect$. We write
$\C_\ad$ for the adjoint subcategory and $U_\C$ for the universal grading group,
whose trivial component is $\C_\ad$; it is abelian if $\C$ is braided
\cite[Proposition~2.3, Remark~2.7]{DGNO}.

\subsection{Transfer of property F}\label{ssec:transfer}
Property F passes to direct summands and fusion subcategories. Indeed, if $V$ is
a summand of $X$, naturality of the braiding makes $V^{\otimes n}$ invariant
under the braid action on $X^{\otimes n}$, and restriction gives a quotient of
$G_{X,n}$. It also passes along a dominant braided tensor functor: the braid
image on the image of an object is a quotient of its original braid image, and
every object of the target is a summand of such an image.

If $\C$ has property F, so does $\C\boxtimes\C^\rev$. Reversing the braiding
inverts the braid generators. On $X\boxtimes Y$, the braid image is a quotient
of a subgroup of $G_{X,n}\times G_{Y,n}$, and every object of the product
category is a summand of some $X\boxtimes Y$.

\subsection{Pivotal dimensions}\label{ssec:dims}
Let $\C$ be pivotal. We write $\Tr$ for the right categorical trace,
$d_X=\Tr(\id_X)$, and $\tr$ for the ordinary matrix trace. For $T\in\End(Z)$,
\[
\Tr(T)=\sum_Wd_W\tr\bigl(T|_{\Hom(W,Z)}\bigr).
\]
For $T\in\End(Y\otimes Z)$, the right partial trace $\ptr_Z(T)\in\End(Y)$
satisfies $\Tr\circ\ptr_Z=\Tr$, $\ptr_Z((R\otimes\id_Z)T)=R\,\ptr_Z(T)$ and
$\ptr_Z(\id_{Y_1}\otimes T_1)=\id_{Y_1}\otimes\ptr_Z(T_1)$. If $Y$ is simple,
then $\ptr_Z(T)=(\Tr(T)/d_Y)\id_Y$; here $d_Y\ne0$, since the dimension of a
simple object in a pivotal fusion category is nonzero \cite[Proposition~4.8.4]{EGNO}.

The function $X\mapsto d_X$ is a character of the fusion ring, so
Perron--Frobenius theory gives
\[
|d_X|\le\FP(X).
\]
For simple $X$, the left dimension of $X$ is $d_{X^*}$ \cite[\S2]{La}, and
$d_{X^*}=\overline{d_X}$ \cite[Proposition~2.9]{ENO}. Hence the squared norm of a
simple object $X$ is $|X|^2=d_Xd_{X^*}=|d_X|^2$, and
$\dim(\C)=\sum_X|d_X|^2\le\FP(\C)$ \cite[Proposition~8.21]{ENO}, the sum running
over the simple objects. For simple $X$, the scalars
$\mathrm{Sl}_X=(d_{X^*}/d_X)\,\id_X$ extend to a tensor automorphism $\mathrm{Sl}$
of the identity functor \cite[proof of Proposition~2.13]{La}, the slope. So if $X$
is simple and $\kappa=d_X/|d_X|$, then $\overline{d_Y}/d_Y=\kappa^{-2n}$ for
every simple summand $Y$ of $X^{\otimes n}$, that is,
\begin{equation}\label{eq:phase}
d_Y\in\kappa^n\mathbb R .
\end{equation}

Pivotal structures form a torsor over
$\Aut_\otimes(\id_\C)\cong\Hom(U_\C,\mathbb C^\times)$ \cite[Exercise~4.7.16]{EGNO},
\cite[Proposition~2.3]{DGNO}, and twisting by $\psi$ multiplies $d_X$ by
$\psi(\deg X)$. On $\C_\ad$ the slope is trivial, so $d_Y$ is real there and
does not depend on the pivotal structure. A pseudo-unitary category has a
unique pivotal structure with $d_X=\FP(X)$, and it is spherical
\cite[Proposition~8.23]{ENO}. A weakly integral category is pseudo-unitary
\cite[Proposition~8.24]{ENO}.

\subsection{Galois conjugation}\label{ssec:galois}
A braided fusion category, possibly with a pivotal structure, has a model over
$\Qb$. The fusion category has a model over a number field
\cite[Theorem~2.30 and Remark~2.33]{ENO}, and braidings and pivotal structures are finitely many
solutions of polynomial equations. For $\sigma\in\Gal(\Qb/\mathbb Q)$, the
conjugate $\C^\sigma$ has the same fusion rules and Frobenius--Perron
dimensions. Its braid matrices are obtained by applying $\sigma$ entrywise, so
property F is invariant under Galois conjugation, and centralisers, $\C'$ and
$\C_\ad$ correspond. If $\C$ is pivotal, so is $\C^\sigma$, with dimensions
$\sigma(d_X)$. The squared norm is defined in any fusion category
\cite[\S2]{ENO} and is computed in the $\Qb$-model, so
$\dim(\C^\sigma)=\sigma(\dim\C)$. We do not use the form $|d_X|^2$ here, since
complex conjugation need not commute with $\sigma$.

\subsection{The Frobenius--Perron object}\label{ssec:FP}
Let $\M$ be a pivotal braided fusion category and put
$S_{XY}=\Tr(c_{Y,X}c_{X,Y})$. The matrix $S$ is symmetric, and for simple $X$
the function $\chi_X(Y)=S_{XY}/d_X$ is a character of the fusion ring
\cite[\S2.2, Proposition~2.4]{La}. For a fusion subcategory $\D$ and a simple object
$V$,
\begin{equation}\label{eq:centraliser}
V\in\D'\iff S_{YV}=d_Yd_V\ \text{ for every }Y\in\cO(\D).
\end{equation}
This is \cite[Proposition~C.3]{DGNO}, whose normalised matrix $\tilde s_{YV}$ equals
$S_{YV}/(d_Yd_V)$ by \cite[\S2.2]{La}; no sphericality is needed.

\begin{proposition}\label{prop:FPobject}
Let $\M$ be non-degenerate and pivotal. There is a unique simple object $P$
with $\chi_P=\FP$, and
\begin{align}
S_{XP}&=d_P\FP(X)\quad\text{for every simple }X,\label{eq:FPcol}\\
\dim(\M)&=|d_P|^2\FP(\M).\label{eq:FPdim}
\end{align}
\end{proposition}

\begin{proof}
Non-degeneracy makes the $\chi_X$ distinct, so they exhaust the characters
\cite[\S2.4]{La}, and $\FP$ is one of them \cite[Proposition~3.3.6]{EGNO}. This gives
$P$ and~\eqref{eq:FPcol}. By \cite[Corollary~2.20]{La}, the matrix $S$ becomes
unitary after division by $\sqrt{\dim\M}$ and by a scalar, which has absolute
value one by \cite[Proposition~2.13]{La}. Hence $SS^\dagger=\dim(\M)\,I$, and
the row of $P$ gives
\[
\dim(\M)=\sum_X|S_{PX}|^2=|d_P|^2\sum_X\FP(X)^2=|d_P|^2\FP(\M).\qedhere
\]
\end{proof}

We call $P$ the \emph{Frobenius--Perron object} of $\M$.

\section{The trace argument}\label{sec:trace}

Fix an object $X$ of a pivotal braided fusion category for which every
$G_{X,n}$ is finite, and let $A_n\subseteq\End(X^{\otimes n})$ be the linear
span of $G_{X,n}$. Averaging makes the group action unitary on each
multiplicity space. The resulting adjoint makes $A_n$ a finite-dimensional
$*$-algebra, with
\[
\Bigl(\sum_{h\in G_{X,n}}\alpha_hh\Bigr)^*=\sum_{h\in G_{X,n}}\overline{\alpha_h}\,h^{-1}.
\]
Thus the involution is independent of the averaged inner products. Tensoring on
either side with identity strands preserves it, since it carries braid images
and their inverses to braid images and inverses.

The next lemma obtains a trace bound by placing overlapping copies of a unitary
on increasingly many strands.

\begin{lemma}[Trace bound]\label{lem:trace}
Let $X$ be simple in a pivotal braided fusion category, with $G_{X,n}$ finite
for every $n$. Put $d=d_X$ and $f=\FP(X)$. Every unitary $u\in A_m$, $m\ge2$,
satisfies
\[
|\Tr(u)|\le|d|f^{m-1}.
\]
\end{lemma}

\Needspace{13\baselineskip}
\begin{proof}
Place $k$ copies $u_1,\dots,u_k$ on consecutive blocks of $m$ strands,
successive blocks sharing one strand, so $\ell=k(m-1)+1$.
\begin{center}
\centering
\begin{tikzpicture}[x=0.44cm,y=0.48cm,line width=0.55pt]
  \foreach \x in {0,...,6}{\draw[->] (\x,0)--(\x,5.2);}
  \draw[fill=white] (-0.28,0.65) rectangle (2.28,1.45);
  \node at (1,1.05) {$u_1$};
  \draw[fill=white] (1.72,2.05) rectangle (4.28,2.85);
  \node at (3,2.45) {$u_2$};
  \draw[fill=white] (3.72,3.45) rectangle (6.28,4.25);
  \node at (5,3.85) {$u_3$};
  \node[below] at (0,0) {$1$};
  \node[below] at (6,0) {$\ell$};
\end{tikzpicture}

\vspace{-2pt}
{\footnotesize Adjacent blocks share one strand.}
\end{center}

These insertions preserve the involution, so $U=u_1\cdots u_k$ is unitary.
On each multiplicity
space its ordinary trace has absolute value at most the dimension of that
space. Consequently,
\[
|\Tr(U)|=\Bigl|\sum_Wd_W\tr\bigl(U|_{\Hom(W,X^{\otimes\ell})}\bigr)\Bigr|
\le\sum_WN^W_{X^{\otimes\ell}}|d_W|\le\sum_WN^W_{X^{\otimes\ell}}\FP(W)=f^\ell.
\]
Simplicity of $X$ makes the partial trace of $u$ over its last $m-1$ strands
equal to $(\Tr(u)/d)\id_X$. Only $u_k$ meets the last $m-1$ strands of $U$;
tracing them out and iterating gives
\[
\Tr(U)=\frac{\Tr(u)^k}{d^{k-1}},\qquad |\Tr(u)|^k\le|d|^{k-1}f^{k(m-1)+1}.
\]
Take $k$-th roots and let $k\to\infty$.
\end{proof}

\begin{lemma}[Saturation]\label{lem:saturation}
Let $X$ be simple in a pivotal braided fusion category, with $G_{X,n}$ finite
for every $n$. Put $d=d_X$, $f=\FP(X)$, $c=c_{X,X}$ and $t=\Tr(c^2)/d$. Then
$|t|\le f$, and if $|t|=f$ then $|d|=f$.
\end{lemma}

\begin{proof}
Since $t=\chi_X(X)$, the bound $|t|\le f$ holds. Assume $|t|=f$ and put
$\kappa=d/|d|$. We first bound the total absolute weight of the categorical
trace on $A_3$. Decompose the multiplicity spaces into irreducible
representations of the finite group $G_{X,3}$:
\[
H_W=\Hom(W,X^{\otimes3})\cong\bigoplus_rV_r^{\oplus m_{Wr}},\qquad
A_3\cong\bigoplus_r\End(V_r),
\]
retaining only the constituents that occur. For $T\in A_3$,
\[
\Tr(T)=\sum_rw_r\tr(T|_{V_r}),\qquad w_r=\sum_Wm_{Wr}d_W .
\]
By~\eqref{eq:phase}, $w_r=\kappa^3v_r$ with $v_r$ real. Let $e\in A_3$ be the
central idempotent onto the blocks with $w_r\ne0$, and let $J\in A_3$ act on
$V_r$ by $\sgn(v_r)$ if $w_r\ne0$ and by $1$ otherwise. This is a unitary
involution, so Lemma~\ref{lem:trace} gives
\begin{equation}\label{eq:J}
\sum_r|w_r|\dim V_r=|\Tr(J)|\le|d|f^2.
\end{equation}

\noindent\begin{minipage}[t]{0.6\textwidth}
\setlength{\parindent}{1em}\indent
We now exhibit a braid whose trace attains this bound. Let $a,b$ be the
images of $\sigma_1,\sigma_2$. Since $X$ is simple, $\ptr_X(c^2)=t\,\id_X$, so
tracing over the third strand of $a^2b^2$ gives $t\,c^2$, and
\[
\Tr(a^2b^2)=t\Tr(c^2)=dt^2,\qquad |\Tr(a^2b^2)|=|d|f^2.
\]
\end{minipage}\hfill
\begin{minipage}[t]{0.37\textwidth}
\vspace{-4pt}
\centering
\begin{tikzpicture}[scale=0.42,line width=0.6pt]
\crossing{1}{0}\crossing{1}{1}\crossing{0}{2}\crossing{0}{3}
\draw (0,0)--(0,2);\draw[->] (0,4)--(0,4.6);
\draw[->] (1,4)--(1,4.6);
\draw (2,2)--(2,4) .. controls (2,5.2) and (3.3,5.2) .. (3.3,4)
  -- (3.3,0) .. controls (3.3,-1.2) and (2,-1.2) .. (2,0);
\node at (4.6,2) {$=\;t$};
\begin{scope}[xshift=5.9cm,yshift=0.5cm]
\crossing{0}{0.5}\crossing{0}{1.5}
\draw (0,0)--(0,0.5) (1,0)--(1,0.5);
\draw[->] (0,2.5)--(0,3.1);\draw[->] (1,2.5)--(1,3.1);
\end{scope}
\end{tikzpicture}

\vspace{1pt}
{\footnotesize Partial trace over strand 3.}
\end{minipage}

\smallskip
Since $g=a^2b^2$ is unitary,
\[
|d|f^2=\Bigl|\sum_rw_r\tr(g|_{V_r})\Bigr|\le\sum_r|w_r|\,|\tr(g|_{V_r})|\le
\sum_r|w_r|\dim V_r\le|d|f^2.
\]
Equality holds throughout. Where $w_r\ne0$, the unitary matrix $g|_{V_r}$ has
trace of absolute value $\dim V_r$, and hence is scalar. Hence $ge$ is central
in $A_3$. Since $e$ is central and $b^2$ commutes with $b$,
\[
0=[a^2b^2e,\,b]=[a^2e,\,b]\,b^2,
\]
so $a^2e$ commutes with $b$. It also commutes with $a$, so it is central in
$A_3$. The braid relation gives $b=(ab)\,a\,(ab)^{-1}$, and conjugating the
central element $a^2e$ by $ab$ yields $b^2e=a^2e$. Thus $a^2b^{-2}e=e$, so
$a^2b^{-2}$ acts trivially on every block of nonzero weight, and
\begin{equation}\label{eq:d3}
\Tr(a^2b^{-2})=\sum_rw_r\dim V_r=\Tr(\id_{X^{\otimes3}})=d^3=\kappa^3|d|^3.
\end{equation}

We compute the same trace directly. Since $c^2$ has finite order, its inverse
has the conjugate trace on each $L_W=\Hom(W,X\otimes X)$, and
$d_W=\kappa^4\overline{d_W}$ by~\eqref{eq:phase}. Hence
$\Tr(c^{-2})=\kappa^4\,\overline{\Tr(c^2)}=\kappa^4\bar t\bar d$. Tracing over
the third strand of $a^2b^{-2}$ gives $(\Tr(c^{-2})/d)\,c^2$, so
\[
\Tr(a^2b^{-2})=\frac{\Tr(c^{-2})\Tr(c^2)}{d}=\kappa^4|t|^2\bar d=\kappa^3|d|f^2 .
\]
Comparing with~\eqref{eq:d3} gives $|d|=f$.
\end{proof}

\begin{remark}
The phase $\kappa$ is needed because pivotal dimensions may be non-real. In a
spherical category $\kappa=\pm1$, and the weights $w_r$ are real. If the
chosen pivotal structure has positive dimensions, for example the canonical
pivotal structure on a pseudo-unitary category, every weight is nonzero and
$e=1$.
\end{remark}

\section{Proof of the main result}\label{sec:main}

\begin{proposition}\label{prop:nondeg}
Let $\M$ be a non-degenerate pivotal braided fusion category with
Frobenius--Perron object $P$. If $G_{P,n}$ is finite for every $n$, then $\M$ is
pseudo-unitary and $P$ is invertible.
\end{proposition}

\begin{proof}
By~\eqref{eq:FPcol}, $t=S_{PP}/d_P=\FP(P)$, so Lemma~\ref{lem:saturation}
gives $|d_P|=\FP(P)\ge1$. By~\eqref{eq:FPdim},
\[
\dim(\M)=|d_P|^2\FP(\M)\ge\FP(\M)\ge\dim(\M).
\]
Hence $\M$ is pseudo-unitary and $\FP(P)=|d_P|=1$.
\end{proof}

\begin{proof}[Proof of Theorem~\ref{thm:main}]
By Proposition~\ref{prop:nondeg}, $\C$ is pseudo-unitary. Let $\sigma$ be a
Galois automorphism. The conjugate $\C^\sigma$ is again non-degenerate and
pivotal with property F (\S\ref{ssec:galois}), hence pseudo-unitary. Their
common fusion rules give
\[
\sigma(\dim\C)=\dim(\C^\sigma)=\FP(\C^\sigma)=\FP(\C)=\dim(\C).
\]
Thus every Galois automorphism fixes $\dim(\C)$. As $\dim(\C)$ is an algebraic
integer \cite[Theorem~2.3, Remark~2.5]{ENO}, $\FP(\C)=\dim(\C)\in\Z$.
\end{proof}

\section{Degenerate categories and consequences}\label{sec:degenerate}

To treat a degenerate category $\C$, we place the Frobenius--Perron object of
its Drinfeld centre in a subcategory with property F. As in the
non-degenerate case, the dimensions of $\C$ need not be real; we need only a
sign condition on $\C'\cap\C_\ad$, where dimensions are real and independent
of the pivotal structure (\S\ref{ssec:dims}).

Let $\C$ be pivotal, and write an object of the centre $\Zc(\C)$ as
$(V,\gamma)$ with $\gamma_X:V\otimes X\to X\otimes V$. By naturality and
monoidality of the pivotal structure $a$, each $a_V$ is a morphism in
$\Zc(\C)$, so $\Zc(\C)$ is pivotal and the forgetful functor preserves traces
\cite[Exercise~7.13.6]{EGNO}.
The centre is non-degenerate \cite[Corollary~3.9]{DGNO}, with
$\dim\Zc(\C)=\dim(\C)^2$ \cite[Theorem~2.15]{ENO} and
$\FP(\Zc(\C))=\FP(\C)^2$ \cite[Theorem~7.16.6]{EGNO}. Let
\[
\iota_+:\C\longrightarrow\Zc(\C),\ X\longmapsto(X,c_{X,-}),\qquad
\iota_-:\C^\rev\longrightarrow\Zc(\C),\ X\longmapsto(X,c^{-1}_{-,X}),
\]
and denote their images by $\C_+$ and $\C_-$. These are fully faithful braided
tensor functors \cite[Proposition~8.6.1]{EGNO}, and $\iota_+$ preserves both
dimensions and Frobenius--Perron dimensions \cite[Proposition~4.5.7]{EGNO}.

\begin{lemma}\label{lem:centre}
Let $\B=\C_+\vee\C_-$ and $\E=\iota_+(\C')$. Then $\B'=\E$ and $\E'=\B$. If $\C$
has property F, so does $\B$.
\end{lemma}

\begin{proof}
The canonical images are mutual centralisers: the double braiding of
$\iota_+(X)$ with $(V,\gamma)$ is $\gamma_X\circ c_{X,V}$, which is the identity
for every $X$ precisely when $(V,\gamma)=\iota_-(V)$, and similarly for $\C_-$.
Consequently $\B'=\C_+\cap\C_-=\iota_+(\C')$, since $\iota_+(X)\cong\iota_-(Y)$
forces $X\cong Y$ and $c_{X,Z}=c^{-1}_{Z,X}$ for all $Z$. The double centraliser
theorem gives $\E'=\B''=\B$ \cite[Corollary~3.11]{DGNO}. Finally, the two canonical
images centralise each other, so $X\boxtimes Y\mapsto\iota_+(X)\otimes\iota_-(Y)$
defines a braided functor $\C\boxtimes\C^\rev\to\B$ \cite[Proposition~8.6.1]{EGNO}. It
is dominant because these images generate $\B$ and commute up to isomorphism.
The transfer properties in \S\ref{ssec:transfer} therefore give property F for
$\B$.
\end{proof}

\begin{lemma}\label{lem:normalise}
Let $\C$ be a pivotal braided fusion category with $d_Y>0$ for every simple
$Y\in\C'\cap\C_\ad$. Then $\C$ has a pivotal structure with $d_Y=\FP(Y)$ for
every simple $Y\in\C'$.
\end{lemma}

\begin{proof}
The symmetric category $\C'$ is integral \cite[Corollary~2.46]{DGNO}, hence
pseudo-unitary, with a pivotal structure of dimensions $\FP$
\cite[Propositions~8.23, 8.24]{ENO}. The restricted pivotal structure differs from it
by a character $\eta$ of $U_{\C'}$, so $d_Y=\eta(\deg Y)\FP(Y)$ for simple
$Y\in\C'$. Inclusion induces a homomorphism $\phi:U_{\C'}\to U_\C$, and since
the universal grading is faithful with trivial component $\C_\ad$, the kernel of
$\phi$ consists of the degrees of the simple objects of $\C'\cap\C_\ad$. For
$g\in\ker\phi$, choose a simple $Y\in\C'$ of degree $g$; then
$Y\in\C'\cap\C_\ad$ and $\eta(g)=d_Y/\FP(Y)>0$. Since $U_{\C'}$ is finite, the
values of $\eta$ are roots of unity, so positivity makes $\eta$ trivial on
$\ker\phi$. Thus $\eta$ factors through $\phi$ and extends to
a character $\psi$ of the finite abelian group $U_\C$. Twisting the pivotal
structure by $\psi^{-1}$ gives $d_Y=\FP(Y)$ on $\C'$ (\S\ref{ssec:dims}).
\end{proof}

\begin{theorem}\label{thm:degenerate}
Let $\C$ be a pivotal braided fusion category with property F such that
$d_Y>0$ for every simple object $Y$ of $\C'\cap\C_\ad$. Then $\C$ and all its
Galois conjugates are pseudo-unitary, and $\dim(\C)=\FP(\C)\in\Z$.
\end{theorem}

\begin{proof}
By Lemma~\ref{lem:normalise} we may assume that $d_Y=\FP(Y)$ for simple
$Y\in\C'$. Set $\M=\Zc(\C)$, let $\B$ and $\E$ be as in Lemma~\ref{lem:centre},
and let $P$ be the Frobenius--Perron object of $\M$. For every simple $Y\in\E$
we have $d_Y=\FP(Y)$, and~\eqref{eq:FPcol} gives
\[
S_{YP}=d_P\FP(Y)=d_Yd_P .
\]
The centraliser criterion~\eqref{eq:centraliser} gives $P\in\E'=\B$. Since $\B$
has property F, all braid images of $P$ are finite, and
Proposition~\ref{prop:nondeg} shows that $\M$ is pseudo-unitary. Thus
$\dim(\C)^2=\FP(\C)^2$, and as both are positive, $\C$ is pseudo-unitary.

Let $\sigma$ be a Galois automorphism. The conjugate $\C^\sigma$ is pivotal
with property F, and its dimensions on $(\C^\sigma)'=(\C')^\sigma$ are
$\sigma(\FP(Y))=\FP(Y)$, since $\C'$ is integral. So $\C^\sigma$ satisfies the
hypothesis and is pseudo-unitary. As in the proof of Theorem~\ref{thm:main},
$\dim(\C)=\FP(\C)\in\Z$.
\end{proof}

The hypothesis holds when $\C$ is non-degenerate, and more generally when
$\C'\cap\C_\ad=\Vect$. It also holds, after conjugation, when $\C$ has a
pseudo-unitary Galois conjugate, and then no pivotal structure need be given.

\begin{corollary}\label{cor:pu}
Let $\C$ be a braided fusion category with a pseudo-unitary Galois conjugate.
If $\C$ has property F, then $\C$ is pseudo-unitary and weakly integral.
\end{corollary}

\begin{proof}
A pseudo-unitary conjugate $\C^\sigma$ with its canonical spherical structure
has $d_Y=\FP(Y)>0$ for every simple $Y$, and it has property F. By
Theorem~\ref{thm:degenerate} it is weakly integral, hence so is $\C$, which is
then pseudo-unitary \cite[Proposition~8.24]{ENO}.
\end{proof}

The theorem also has an objectwise consequence. To apply it to the fusion
subcategory $\langle X\rangle$ generated by an object $X$, we must establish
property F for every object of that subcategory. The following lemma shows that
finiteness of all $G_{X,n}$ suffices.

\begin{lemma}[Closure]\label{lem:closure}
Let $X$ be an object of a braided fusion category. If $G_{X,n}$ is finite for
every $n$, then $\langle X\rangle$ has property F.
\end{lemma}

\begin{proof}
The case $X=0$ is immediate. Otherwise, summands of finite sums of positive
powers of $X$ form a nonzero full abelian subcategory closed under tensor
products and subquotients. By \cite[Corollary~4.11.4]{EGNO}, it contains the unit and
is rigid, so it equals $\langle X\rangle$. Thus every object of $\langle
X\rangle$ is a summand of some $Z=\bigoplus_iX^{\otimes k_i}$, with $k_i\ge1$.
For each $n\ge2$, the pure braid group $P_n$ preserves each summand
$X^{\otimes k_{i_1}}\otimes\cdots\otimes X^{\otimes k_{i_n}}$ of $Z^{\otimes n}$.
By cabling, its action there has image in the finite group
$G_{X,k_{i_1}+\cdots+k_{i_n}}$. Thus its image on $Z^{\otimes n}$ is finite,
being contained in a finite product of finite groups. Since $[B_n:P_n]=n!$, the
full braid image is finite as well \cite[Lemma~2.1]{ERW}. Naturality passes
finiteness to every summand of $Z$, as in \S\ref{ssec:transfer}.
\end{proof}

Consequently, the dimension of a particular simple object can force its own
braid images to be infinite.

\begin{corollary}\label{cor:object}
Let $X$ be a simple object of a braided fusion category, and suppose that
$\langle X\rangle$ is non-degenerate and pivotal, or has a pseudo-unitary
Galois conjugate. If $\FP(X)^2\notin\Z$, then $G_{X,n}$ is infinite for every
sufficiently large $n$.
\end{corollary}

\begin{proof}
If every $G_{X,n}$ were finite, Lemma~\ref{lem:closure} together with
Theorem~\ref{thm:main} or Corollary~\ref{cor:pu} would make $\langle X\rangle$
weakly integral. Then $\FP(X)^2\in\Z$ by \cite[Proposition~8.27]{ENO}, a
contradiction. Hence $G_{X,n_0}$ is infinite for some $n_0$. Tensoring with a
nonzero object is faithful, so tensoring with identity strands embeds
$G_{X,n_0}$ in $G_{X,n}$ for every $n\ge n_0$.
\end{proof}

In particular, the conclusion holds in every unitary braided fusion category,
since its fusion subcategories are unitary and hence pseudo-unitary. The
simplicity hypothesis is essential: direct sums in a weakly integral category
need not have integral squared Frobenius--Perron dimension.

\section*{Acknowledgements}
We thank Clement Delcamp for many illuminating discussions. We are grateful to
\href{https://projectnumina.ai}{Project Numina} and XTX Markets for
computing support, and to Project Numina and Zekai Zhu for collaboration
on related topics. JBD is supported by a Four-Year
Fellowship (4YF) from the University of British Columbia. CC and RT
are supported by the Simons Collaboration on Global Categorical
Symmetries and the US Department of Energy Grant 5-29073.

\begingroup
\linespread{1}\selectfont

\endgroup

\end{document}